\documentclass[12pt, twoside]{article}
 \usepackage{amsfonts, amssymb, amsmath, amsthm, latexsym, amscd}
\newcommand{\braket}[2]{\langle #1,#2 \rangle}

\newcommand{\dDel}{\bigtriangledown}

\newcommand{\Ga}{\Gamma}
\newcommand{\ga}{\gamma}
\newcommand{\Om}{\Omega}
\newcommand{\del}{\delta}
\newcommand{\Del}{\Delta}
\newcommand{\La}{\Lambda}

\newcommand{\tht}{\theta}
\newcommand{\eps}{\epsilon}

\newcommand{\f}{\frac}

\newcommand{\lo}{\longrightarrow}

\newcommand{\BR}{\mathbb{R}}
\newcommand{\al}{\alpha}

\newcommand{\pa}{\partial}

\theoremstyle{plain}
\newtheorem{Th}{Theorem}[section]
\newtheorem{Lem}[Th]{Lemma}
\newtheorem{Prop}[Th]{Proposition}

\theoremstyle{definition}

\newtheorem{Rem}[Th]{Remark}
\begin{document}
\addtolength{\textheight}{0 cm} \addtolength{\hoffset}{0 cm}
\addtolength{\textwidth}{0 cm} \addtolength{\voffset}{0 cm}

\title{Existence of soliton solutions for  a quasilinear Schr\"{o}dinger
equation involving critical exponent in ${\BR}^N$}

\author{Abbas Moameni  \textsc \bf \footnote{Research is supported by a
Postdoctoral Fellowship at the University of British Columbia.}}
\maketitle
\begin{center}
{\small Department of Mathematics \\
\small  University of British Columbia\\
\small Vancouver, B.C., Canada \\
\small {\tt moameni@math.ubc.ca }}
\end{center}

\begin{abstract}
Mountain pass in a suitable  Orlicz space is employed to prove the
existence of soliton solutions for a quasilinear Schr\"{o}dinger
equation involving critical exponent in ${\BR}^N$. These equations
contain strongly singular nonlinearities which include derivatives
of the second order. Such equations have been studied as models of
several physical phenomena. The nonlinearity here corresponds to the
superfluid film equation in plasma physics.
\end{abstract}
\section{Introduction}

We study  the existence of  standing wave solutions for quasilinear
Schr\"{o}dinger equations of the form
\begin{equation}
i\pa_t z=-\epsilon \Del z+W(x)z-l(|z|^2)z-k\epsilon\Del
h(|z|^2)h'(|z|^2)z, \quad x\in {\BR}^N,  N >2,
\end{equation}
where $W(x)$  is a given potential, $k$ is a real constant and $l$
and $h$ are real functions that are essentially pure power forms.
%
%This equation very recently has been considered by the authors in
% ,....  The semilinear case corresponding to $k=0$ has been studied
%extensively in recent years (e.g., [3], [9], [24]).
Quasilinear equations of the form (1) have been established  in
several areas of physics corresponding to various types of $h$.
 The  superfluid film equation in
plasma physics has this structure   for $h(s)=s$, ( Kurihura in
[7]). In the case $h(s)=(1+s)^{1/2}$, equation (1) models the
self-channeling of a high-power ultra short laser in matter, see
 [19]. Equation (1) also
appears in  fluid mechanics [7,8], in the theory of Heidelberg
ferromagnetism and magnus [9], in dissipative quantum mechanics [6]
and in condensed matter theory [13]. We consider the case  $h(s)=s$,
$l(s)=\mu s^{\frac{p-1}{2}}$ and $k>0$. Setting
$z(t,x)=\exp(-iFt)u(x)$ one obtains a corresponding equation of
elliptic type which has the formal variational structure:
\begin{eqnarray}
-\epsilon \Delta u+V(x)u-\epsilon k(\Del(|u|^{2}))u=\mu  |u|^{p-1}u,
\quad u>0, x \in
  {\BR}^N,
\end{eqnarray}
where $V(x)=W(x)-F$ is the new potential function.\\

 Note that
$p+1=22^*= \frac{4N}{N-2}$
 behaves like a critical exponent for the above equation [12, Remark
 3.13]. For the
 subcritical case $p+1 <22^*$ the existence of solutions for problem (2) was studied in [10, 11, 12, 14, 15,
 16] and it was left open for the critical exponent case $p+1=22^*$ [12; Remark
 3.13].
   In the present paper,   the
existence of solutions is proved  for $p+1=22^*$ whenever the
potential function $V(x)$ satisfies some geometry conditions. It is
well-known that  for the semilinear case $(k=0)$, $p+1=2^*$ is the
critical exponent. In this case there are many results about the
existence of solutions for the subcritical  and the critical
exponent (e.g. [1, 4, 5, 17, 20]).\\

 %For the subcritical case  equation (2), very recently, has been considered in
%[10, 11, 12, 14 , 15, 16] but the semilinear case corresponding to
%$k=0$ has been studied more extensively (e.g., [1, 4, 17, 20]).
 In the case $k >0$, for a family of parameter $\mu$,  the
existence of a nonnegative solution is proved  for $N=1$ by
Poppenberg, Schmitt and Wang in [16] and for $N\geq 2$ by  Liu and
Wang  in [11].
%for $ 2<p+1 <2\alpha 2^{*}$ in the compact case and for
%$ 4\alpha \leq p+1 <2\alpha 2^{*}$ for the locally compact case
   In [12] Liu and Wang improved
these results for any $\mu >0$ by using a change of variables and
treating the new problem in an Orlicz space. The author in [14],
using the idea of the fibrering method, studied this problem in
connection with the corresponding eigenvalue problem for the
laplacian $ -\Del u=V(x) u$ and  proved the existence of multiple
solutions for  problem (2).
 It is  established in  [10],
the existence of both one-sign and nodal ground states of soliton
type solutions by the Nehari method. They also established some
regularity of the positive solutions.

In this paper, we assume  that the potential function $V$ is radial,
that is  $V(x)=V(|x|),$ and  satisfies the following conditions:

There exist $0<R_1<r_1<r_2<R_2$ and $\alpha>0$ such that
\begin{align*}
V(x) =0, & \ \ \ \ \  \forall x\in
\Omega : =\left\{ x\in {\mathbb{R}}^N : r_1 < |x| < r_2 \right\},\quad \quad \quad  \quad (A_1) \\
V(x) \geq \alpha, & \ \ \ \ \  \forall x\in\Lambda^c,  \quad \quad
\quad \quad \quad \quad \quad \quad \quad \quad \quad \quad \quad
\quad \quad \quad  (A_2)
\end{align*}
where $\Lambda=\left\{ x \in {\BR}^N : R_1< |x| < R_2 \right\}.$

Here is our main Theorem.
\begin{Th}
 There exists $\eps_0>0$, such that for all $\eps\in
(0,\eps_0)$ problem $(2)$ has a nonnegative solution $u_\eps\in
H^1_r({\BR}^N)$ with $u^2_\eps\in H^1_r({\BR}^N)$ and
$$u_\eps(x)\lo 0 \quad as\quad |x|\lo+\infty.$$
\end{Th}
This paper is organized as follows. In Section 2, we reformulate
this problem in an appropriate Orlicz space. In Section 3, we prove
the existence of a solution for   a special deformation of  problem
(2). Theorem 1.1 is proved in Section 4.
\section{Reformulation of the problem and preliminaries}
Denote by $H_{r}^1 ({\BR}^{N})$ the space of radially symmetric
functions in
\begin{eqnarray*}
H^{1,2} ({\BR}^N)= \left \{ u \in L^{2}({\BR}^N) : \dDel u \in
L^{2}({\BR}^N) \right \}.
\end{eqnarray*}
Denote by $C_{0,r}^{\infty} ({\BR}^N)$  the space of radially symmetric
functions in $C_{0}^{\infty} ({\BR}^N)$ and  by $D^{1,2} ({\BR}^N)$ the following space,
\begin{eqnarray*}
D^{1,2} ({\BR}^N):=\left \{ u \in L^{2^*} ({\BR}^N) :   \dDel u \in L^{2}({\BR}^N) \right \}.
\end{eqnarray*}
Without loss   of generality, one can assume $k=1$ in problem (2).
 We formally
formulate problem (2) in a variational structure as follows
$$J_\eps(u)=\f{\eps}{2} \int_{{\BR}^N}(1+u^2)| \dDel u|^2dx+\f{1}{2}
\int_{{\BR}^N}V(x)u^2dx-\f{1}{22^*} \int_{{\BR}^N} |u|^{22^*}dx$$ on
the space
$$X=\{u\in H_r^{1,2}({{\BR}}^N): \int_{{{\BR}}^N} V(x) u^2dx<\infty\},$$
which is equipped with the following norm,
$$\|u\|_X=\left \{ \int_{{\BR}^N}| \dDel u|^2dx+ \int_{{{\BR}}^N} V(x) u^2dx\right \}^{\frac{1}{2}}.$$
 Liu and Wang in [12] for the subcritical case, i.e.
$$J_\eps(u)=\f{\eps}{2} \int_{{\BR}^N}(1+u^2)| \dDel u|^2dx+\f{1}{2}
\int_{{\BR}^N}V(x)u^2dx-\f{1}{p+1} \int_{{\BR}^N} |u|^{p+1}dx,\quad
\quad (2< p+1 < 22^*),$$  by making a change of variables treated
this problem in an  Orlicz space. Following their work, we consider
this problem for the critical exponent case $(p+1 = 22^*)$ in the
same Orlicz space. To  convince the reader we briefly recall some of
their notations and results that are useful in the sequel.

First, we make a change of variables as follows,
$$dv=\sqrt{1+u^2}du, \quad v=h(u)=\f{1}{2}u\sqrt{1+u^2}+\f{1}{2} \ln
(u+\sqrt{1+u^2})$$ Since $h$ is strictly monotone it has a
well-defined inverse function: $u=f(v)$. Note that
$$h(u)\sim \begin{cases}
u, & |u|\ll 1\\
\f{1}{2}u|u|, & |u|\gg 1, \quad h'(u)=\sqrt{1+u^2},
\end{cases}$$
and
$$f(v)\sim \begin{cases}
v & |v|\ll 1 \\
\sqrt{\f{2}{|v|}}v, & |v|\gg 1, \quad
f'(v)=\f{1}{h'(u)}=\f{1}{\sqrt{1+u^2}}= \f{1}{\sqrt{1+f^2(v)}}.
\end{cases}$$
Also, for some $C_0>0$ it holds
$$G(v)=f(v)^2\sim \begin{cases} v^2& |v|\ll 1,\\
2|v|& |v|\gg 1, \quad G(2v)\leq C_0 G(v),
\end{cases}$$
$G(v)$ is convex,
$G'(v)=\f{2f(v)}{\sqrt{1+f(v)^2}}$, $G''(v)=\f{2}{(1+f(v)^2)^2}>0$.

Using this change of variable, we can rewrite the functional
$J_\eps(u)$ as
$$\bar {J}_\eps(v)=\f{\eps}{2} \int_{{\BR}^N} |\dDel v|^2dx+ \f{1}{2} \int_{{\BR}^N}V(x)
f(v)^2 dx-\f{1}{22^*} \int_{{\BR}^N} |f(v)|^{22^*}dx.$$ $\bar
{J}_\eps$ is defined on the space
$$H^1_G ({\BR}^N)=\{v|   v(x)=v(|x|),   \int_{{\BR}^N}|\dDel v|^2dx <\infty, \int_{{\BR}^N}V(x)
G(v)dx<\infty\}.$$ We introduced the Orlicz space (e.g.[18])
$$E_G ({\BR}^N)=\{v| \int_{{\BR}^N}V(x) G(v)dx<\infty\}$$
equipped with the  norm
$$|v|_G=\inf_{\zeta>0} \zeta(1+\int_{{\BR}^N}(V(x) G(\zeta^{-1}v(x))dx),$$
and define the norm of $H^1_G ({\BR}^N)$ by
$$\|v\|_{H^1_G ({\BR}^N)}= |\dDel v|_{L^2({\BR}^N)}+|v|_G$$
Here are  some related facts.
\begin{Prop}
 \begin{enumerate}
 \item[{ (i)}]
 $E_G({\BR}^N)$ is a Banach space.
 \item[{(ii)}]  If $v_n\lo v$ in $E_G({\BR}^N)$, then $\int_{{\BR}^N}V(x)|
G(v_n)-G(v)| dx\lo 0$ and  $\int_{{\BR}^N}V(x)| f(v_n)-f(v)|^2dx\lo
0$.
 \item[{(iii)}]  If $v_n\lo v$ a.e. and $\int_{{\BR}^N}V(x) G(v_n)dx\lo
\int_{{\BR}^N}V(x)G(v)dx$, then $v_n\lo v$ in $E_G({\BR}^N)$.
 \item[{(iv)}]  The dual space $E^*_G({\BR}^N)=L^\infty\cap L_V^2=\{w| w\in L^\infty,
\int_{{\BR}^N}V(x)w^2dx<\infty\}$.
 \item[{ (v)}]  If $v\in E_G({\BR}^N)$, then $w=G'(v)=2f(v)f'(v)\in E^*_G({\BR}^N)$, and
$|w|_{E^*_G}=\sup_{|\phi|_{G\leq 1}}(w,\phi)\leq
C_1(1+\int_{{\BR}^N} V(x)G(v)dx)$, where $C_1$ is a constant
independent of $v$.
 \item[{(vi)}] For $N>2$ the map:$v\lo f(v)$ from $H^1_G({\BR}^N)$ into
$L^q({\BR}^N)$ is continuous for $2\leq q\leq 22^*$ and is compact
for $2< q <22^*.$
 \item[{(vii)}] Suppose $B_k$ is the ball with center at  the coordinate origin  and radius  $k>0.$   Let $r<s$ and $\textit{Q}=B_s \backslash B_r .$   The map:$v\lo f(v)$ from $H^1_G({\BR}^N)$ into
$L^q(\textit{Q})$  is compact for $q\geq 2.$
 \end{enumerate}
\end{Prop}
\paragraph{\bf Proof.} See Propositions~(2.1) and (2.2) in [12] for the proof of parts $(i)$ to $(vi)$. We
prove part $ \it {(vii)}$. Set $u=f(v)$.  It is easy to check that
$\|u\|_X \leq \|v\|_{H^1_G({\BR}^N)}.$
 It is standard that the embedding from  $H^1_r({\BR}^N)$  to $L^q(\textit{Q})$ is compact  for $q \geq 2$ (e.g. [1]).
  Hence, we obtain the map:$v\lo f(v)$ from $H^1_G({\BR}^N)$ into
$L^q(\textit{Q})$  is
compact for $q\geq 2$.   $\square$\\

 Hence forth, $\int, H^1, H^1_r,  H^1_G, E_G, L^t$  and $\|\cdot\|$
stand for $\int_{{\BR}^N}$, $H^{1,2}({\BR}^N)$,  $ H^1_r({\BR}^N)$,
$ H^1_G({\BR}^N)$, $ E_G({\BR}^N)$, $L^t({\BR}^N)$ and
$\|\cdot\|_{H^1_G({\BR}^N)}$ respectively. In the following we use
$C$ to denote any constant that is independent of the sequences
considered.

\section{Auxiliary Problem}
In this section, we shall show some results needed  to prove Theorem
1.1. Indeed,  we first consider a special deformation $\bar{H}_\eps$
(See (3) in the following) of $\bar{J}_\eps.$ Then, We show that the
functional $\bar{H}_\eps$ satisfies all the properties of the
Mountain Pass Theorem. Consequently,  $\bar{H}_\eps$ has a critical
point for each $\eps>0.$  We shall use this to prove Theorem 1.1 in
the next section. In fact, we will see that  the functionals
$\bar{J}_\eps$ and $\bar{H}_\eps$  will coincide for the small
values of $\eps$. This idea was explored by
Del Pino and Felmer [5].\\

To do this, we shall consider constants $\tht,k$ and $\beta$
satisfying
$$4<\tht<22^*, \quad k>\f{\tht}{\tht-2}, \quad
\beta=(\f{\al}{k})^{\f{1}{2(2^*-1)}}, \quad (\al  \quad  \text{is
introduced in}  \quad  A_2)$$ and functions
\begin{align*}
\ga(s)&=\begin{cases}s^{22^*-1}, & s>0,\\ 0, & s\leq 0 \\
\end{cases}\\
\bar{\ga}(s)&=\begin{cases} \ga(s), & s\leq \beta, \\ (\f{\al}{k})s,
& s>\beta,
\end{cases}\\
w(x,s)&=\chi_\La(x)\ga(s)+(1-\chi_\La(x))\bar{\ga}(s),
\end{align*}
where $\chi_\La$ denotes the characteristic function of the set
$\La$. Set $W(x,t)=\int_0^t w(x,\zeta)d\zeta$.
 It is easily seen that the function $w$
satisfies the following conditions,
\begin{eqnarray*}
& 0\leq \tht W(x,t) \leq  w(x,t)t, \quad \forall x\in\La, t\geq 0,\quad \quad \quad \quad \quad \quad (g_1)\\
& 0\leq 2W(x,t)\leq w(x,t)t\leq \f{1}{k} V(x)t^2,\quad
 \forall x\in {\La}^c, t \in {\BR}.  \quad \quad (g_2)
\end{eqnarray*} Now, we study the existence
of solutions for  the deformed equation, i.e.
$$-\eps\Del u+V(x)u-\eps(\Del(|u|^2))u = w(x,u), \quad x\in {\BR}^N.$$
which correspond to the critical points of
$$H_{\eps}(u)=\f{\eps}{2}\int (1+u^2)|\dDel u|^2+\f{1}{2} \int V(x)u^2-\int W(x,u)dx.$$
As in Section (2), we can rewrite the functional  $H_{\eps} (u)$ as
a new functional  $\bar H_{\eps} (v)$ with $u=f(v)$ as follows,
\begin{equation}
\bar{H}_{\eps}(v)=\f{\eps}{2} \int |\dDel v|^2 dx+\f{1}{2} \int
V(x)f(v)^2dx-\int W(x,f(v))dx.
\end{equation}
 $\bar {H}_{\eps} (v)$
is defined on the Orlicz space $H^G_1.$ To simplify the writing in
this section, we shall assume $\eps=1, H_1=H$ and $\bar {H}_1=\bar
H.$

The following Proposition states some properties of the functional
$\bar {H}.$
 \begin{Prop}
 \begin{enumerate}
 \item[{(i)}]  $\bar{H}$ is well-defined on $H^1_G$.
 \item[{(ii)}]  $\bar{H}$ is continuous in $H^1_G$.
 \item[{(iii)}] $\bar{H}$ is Gauteaux-differentiable in $H^1_G$.
 \end{enumerate}
\end{Prop}
%For $v\in H^1_G$, the
%$G$-derivative $I'(v)$ is a continuous linear functional, and
%$I'(v)$ is continuous in $v$ in the strong-weak topology, that is,
%if $v_n\lo v$ strongly in $H^1_G$ then $I'(v_n)\lo I'(v)$ weakly.
\paragraph{\bf Proof.} The proof is similar to the proof of
Proposition~(2.3) in [12] by some obvious changes.$\square$\\

Here is the main result in this section.
\begin{Th}
 $ \bar H$ has a critical point  in $H^G_1$, that is,
there exists $0 \neq v\in H^G_1$ such that
$$\int \dDel
v.\dDel \phi dx+\int V(x)f(v)f'(v)\phi dx-\int w(x,f(v))f'(v)\phi
dx=0,$$ for every $\phi\in H^G_1$.
\end{Th}

We use the Mountain Pass Theorem (see [2], [17]) to prove Theorem
3.2.  First, let us define the Mountain Pass value,
$$C_0:=\inf_{\ga\in\Ga} \sup_{t\in[0,1]} \bar{H}(\ga(t)),$$
where
$$\Ga=\{\ga\in C([0,1], H^1_G) | \ga(0)=0, \bar {H}(\ga(1))\leq 0, \ga(1)\neq 0\}.$$
The following Lemmas are crucial for the proof of Theorem 3.2.
\begin{Lem}
 The functional $\bar{H}$ satisfies the Mountain Pass
Geometry. \end{Lem}
\paragraph{\bf Proof.}
 We need to show that there exists $0\neq v \in~H^1_G$ such that
$\bar{H}(v)\leq 0$.   Let $e\in C_{0,r}^\infty({\BR}^N)$ with
$e\not\equiv 0$ and  supp$( e)\subset\Omega$. It is easy to see
that ${H}(te)\leq 0 $ for the large values of $t.$   Consequently
  $\bar {H} (v) <0$ where $v=h(te).$ $\square$
\begin{Lem} $C_0$\; is positive.
\end{Lem}
\paragraph{\bf Proof.} Set
\begin{equation*}
S_\rho:=\{v\in H^1_G| \int |\dDel v|^2dx+\int V(x) f(v)^2dx=\rho^2\}.
\end{equation*}
For $z=f(v)|f(v)|$ with $v\in S_\rho$, we have
$$\int |\dDel z|^2dx= \int \f{4f^2(v)}{1+f^2(v)} |\dDel v|^2dx\leq 4\rho^2$$
from which, we obtain
\begin{align}
\int |f(v)|^{22^*} dx = \int |z|^{2^*} dx & \leq C(\int |\dDel z|^2
dx)^{\f{2^*}{2}}\notag\\
& \leq C\rho^{2^*}.
\end{align}
Also,  it follows from $(g_1)$  and $(g_2)$ that
\begin{align}
\int W(x,f(v))dx &=\int_\La W(x,f(v))dx+\int_{\La^c} W(x,f(v))dx \notag\\
&\leq \f{1}{\tht} \int |f(v)|^{22^*} dx+ \f{1}{2k}\int V(x)f(v)^2 dx.
\end{align}
Considering  (4), (5) and the fact that $v\in S_\rho$, we obtain
\begin{align*}
\bar{H}(v)&=\f{1}{2}\int |\dDel v|^2dx+\f{1}{2} \int V(x)f(v)^2 dx-\int
W(x,f(v))dx\\
&\geq \f{1}{2}\rho^2-C\rho^{2^*}-\f{1}{2k}
\rho^2=(\f{1}{2}-\f{1}{2k})\rho^2-C\rho^{2^*} \geq \f{k-1}{4k}\rho^2,
\end{align*}
when $0<\rho\leq \rho_0\ll 1$ for some $\rho_0$. Hence, for $v\in
S_\rho$ with $0<\rho\leq \rho_0$ we have
\begin{equation}
\bar{H}(v)\geq \f{k-1}{4k}\rho^2.
\end{equation}
If $\ga(1)=v$ and $\bar {H}(\ga(1))<0$ then it follows from (6) that
$$\int |\dDel v|^2dx+ \int V(x)f(v)^2dx >\rho_0^2,$$
thereby giving
$$\sup_{t\in [0,1]}\bar{H}(\ga(t)) \geq \sup_{\ga(t)\in S_{\rho_0}}
\bar{H}(\ga(t))\geq \f{k-1}{4k}\rho_0^2.$$
Therefore $C_0\geq \f{k-1}{4k}\rho_0^2>0$.$\Box$

The Mountain Pass Theorem guaranties the existence of a $(PS)_{C_0}$
sequence $\{v_n\},$ that is,  $\bar{H}(v_n)\lo C_0$ and
$\bar{H}'(v_n)\lo 0$. The following Lemma states some properties of
this sequence.
\begin{Lem}Suppose $\{v_n\}$ is a $(PS)_{C_0}$ sequence. The following statements hold.
 \begin{enumerate}
 \item[{(i)}] $\{v_n\}$ is bounded in $H^1_G.$
\item[{(ii)}]  For each $\del>0$, there exists $R>4R_2$, ($R_2$ is introduced in $(A_1)$ and $(A_2)$) such that
$$\underset{n\rightarrow+\infty}{\mathrm{limsup}}
\int_{B_R^c} \big ( |\dDel v_n|^2+ V(x)f(v_n)^2 \big )dx<\del.$$
\item[{(iii)}]  If $v_n$ converges weakly to $v $  in $ H^1_G$, then
$$\lim_{n\rightarrow+\infty} \int w(x,f(v_n)) f(v_n)dx=\int w(x,f(v))f(v) dx.$$
\item[{(iv)}] If $v_n\geq 0$ converges  weakly to $v $  in $ H^1_G$,  then for every nonnegative test  function $\phi \in  H^1_G $ we have
$$\lim_{n\rightarrow+\infty} \braket{\bar{H}'(v_n)}{\phi}= \braket{\bar{H}'(v)}{\phi}.$$
\end{enumerate}
\end{Lem}
\paragraph{\bf Proof.} Since $\{v_n\}$ is  a $(PS)_{C_0}$ sequence, we have
\begin{eqnarray}
\bar{H}(v_n)&=&\f{1}{2} \int |\dDel v_n|^2 dx+\f{1}{2} \int V(x)f(v_n)^2
dx- \int W(x,f(v_n))dx\nonumber \\&=&C_0+o(1),
\end{eqnarray}
and
\begin{eqnarray}
\braket{\bar{H}'(v_n)}{\phi}&=&\int \dDel v_n.\dDel \phi dx+\int V(x)f(v_n)f'(v_n)\phi
dx-\int w(x, f(v_n))f'(v_n) \phi dx \nonumber \\
&=&o(\|\phi\|)
\end{eqnarray}
For part $(i),$ pick $\phi=\f{f(v_n)}{f'(v_n)}= \sqrt{1+f(v_n)^2}f(v_n)$ as a test
function. One can easily deduce that $|\phi|_G\leq C|v_n|_G$ and
$$|\dDel \phi|=(1+\f{f(v_n)^2}{1+f(v_n)^2}) |\dDel v_n| \leq 2|\dDel v_n|,$$
which implies  $\|\phi\|\leq C\|v_n\|$. Substituting $\phi$ in (8),
gives
\begin{align}
\braket{\bar{H}'(v_n)}{
 \f{f(v_n)}{f'(v_n)}}&= \int (1+\f{f(v_n)^2}{1+f(v_n)^2}) |\dDel
v_n|^2dx+ \int V(x)f(v_n)^2dx \nonumber \\
&\quad - \int w(x,f(v_n))f(v_n)dx \nonumber \\
&=o(\|v_n\|).
\end{align}
It follows from $(g_1)$ and $(g_2)$ that
\begin{equation}
-\int W(x,f(v_n)) dx+\f{1}{\tht} \int w(x,f(v_n)) f(v_n) dx\geq
\f{1}{k} (\f{1}{\tht}-\f{1}{2}) \int V(x)f(v_n)^2dx
\end{equation}
 Taking into
account (7), (9) and (10), we have
\begin{align*}
C_0+o(1)+o(\|v_n\|) =&\bar{H}(v_n)-\f{1}{\tht}
\braket{\bar{H}'(v_n)}{
\f{f(v_n)}{f'(v_n)}}\\
=&\f{1}{2}\int |\dDel v_n|^2 dx+\f{1}{2} \int V(x)f(v_n)^2dx
-\int
W(x,f(v_n))dx \\
&-\f{1}{\tht} \int (1+\f{f(v_n)^2}{1+f(v_n)^2}) |\dDel v_n|^2 dx -
\f{1}{\tht}
\int V(x)f(v_n)^2dx \\
&+\f{1}{\tht} \int w(x,f(v_n))f(v_n)dx\\
=&\int(\f{1}{2}-\f{1}{\tht}(1+\f{f(v_n)^2}{1+f(v_n)^2})) |\dDel
v_n|^2dx+
(\f{1}{2}-\f{1}{\tht}) \int V(x) f(v_n)^2dx \\
&+\int (\f{1}{\tht} w(x,f(v_n))f(v_n)-W(x,f(v_n)))dx \\
\geq & (\f{1}{2}-\f{2}{\tht}) \int |\dDel v_n|^2
dx+(\f{1}{2}-\f{1}{\tht})(1-\f{1}{k}) \int V(x)f(v_n)^2dx.
\end{align*}
Since, $\f{1}{2}-\f{2}{\tht}>0$ and
$(\f{1}{2}-\f{1}{\tht})(1-\f{1}{k})>0$ it follows from the above
that   $\int |\dDel v_n|^2dx+ \int V(x) f(v_n)^2dx$ is bounded. It
proves part $(i).$

For part $(ii)$,  let $\eta_R\in C^\infty({\BR}^N,{\BR})$ be a
function satisfying $\eta_R=0$ on $B_{\f{R}{2}}$, $\eta_R=1$ on
$B_R^c$ and $|\dDel \eta_R(x)|\leq\f{C}{R}$. It follows from part
$(i)$  that $\{v_n\}$ is bounded. Hence, from (8) we have
$$ \braket{\bar{H}'(v_n)}{ \f{f(v_n)}{f'(v_n)} \eta_R}=o(1),$$
thereby giving
\begin{multline*}
\int (1+\f{f(v_n)^2}{1+f(v_n)^2}) |\dDel v_n|^2\eta_R dx+\int V(x)
f(v_n)^2\eta_R dx\\
+\int \f{f(v_n)}{f'(v_n)} \dDel v_n.\dDel \eta_R dx=\int
w(x,f(v_n))f(v_n)\eta_Rdx+o(1).
\end{multline*}
By $(g_2)$, we get
$$w(x,f(v_n)) f(v_n) \leq \f{V(x)}{k} f(v_n)^2, \quad \forall x\in
B^c_{\f{R}{2}}.$$ Therefore,
\begin{align}
\int (1+\f{f(v_n)^2}{1+f(v_n)^2} )& |\dDel v_n|^2 \eta_R dx +\int
(1-\f{1}{k})V(x) f(v_n)^2\eta_R dx  \nonumber \\
&\leq \f{C}{R} \int \f{|f(v_n)|}{f'(v_n)} |\dDel v_n| dx+o(1) \nonumber  \\
&\leq \f{C}{R} \int |\dDel v_n|^2dx +\f{C}{R} \int
(|f(v_n)|^2+|f(v_n)|^4)dx+o(1).
\end{align}
Also, it follows from part $(vi)$ of Proposition 2.1 that
$\{f(v_n)\}_n$ is a bounded sequence in $L^2({\BR}^N)\cap
L^{22^*}({\BR}^N)$. Hence, $\int (|f(v_n)|^2+ |f(v_n)|^4) dx$ is
bounded. Therefore, it follows from (11) that
$$\underset{n\rightarrow\infty}{\mathrm{limsup}}
\int_{B_R^c} \big  (|\dDel v_n|^2dx+ V(x) f(v_n)^2 \big )dx<\del,
\quad (R>4R_2).$$ It proves part $(ii).$

 For part $(iii),$ note first that from part $(ii)$ of the present   Lemma for each $\del>0$  there exists $R>4R_2$ such that
\begin{equation}
\underset{n\rightarrow\infty}{\mathrm{limsup}}\int_{B_R^c}\big
(|\dDel v_n|^2+V(x) f(v_n)^2 \big )dx<\f{k\del}{4}.
\end{equation}
Since $B_R^c \subseteq \Lambda^c,$ it follows from $(g_2)$ that
$$w(x,f(v_n)) f(v_n)\leq \f{V(x)}{k}f(v_n)^2 \quad \quad \quad \forall x \in B_R^c$$
which together with (12) imply that
\begin{equation}
\underset{n\rightarrow\infty}{\mathrm{limsup}} \int_{B_R^c}
w(x,f(v_n))f(v_n)dx\leq \f{\del}{4},
\end{equation}
 and consequently
\begin{equation*}
\int_{B_R^c} w(x,f(v))f(v)dx\leq \f{\del}{4}.
\end{equation*}
It follows from (13) and the above inequality  that
\begin{multline}
\Big |\int w(x,f(v_n))f(v_n)dx- \int w(x,f(v))f(v)dx \Big |  \\
\leq \f{\del}{2}+\Big | \int_{B_{R_1}} \big [w(x,f(v_n))f(v_n)-w(x,f(v))f(v)\big ]dx \Big |  \\
+\Big | \int_{B_R\backslash B_{R_1}}  \big
[w(x,f(v_n))f(v_n)-w(x,f(v))f(v) \big ]dx \Big |.
\end{multline}
Since $B_{R_1}\subset\La^c$, we have
$$w(x,f(v_n))f(v_n) \leq \f{V(x)}{k} f(v_n)^2, \quad \forall x\in B_{R_1}$$
Then, by the compact theorem embedding and Lebesgue Theorem, we
obtain a subsequence still denoted by $\{v_n\}$, such that
\begin{equation}
\int_{B_{R_1}} w(x,f(v_n))f(v_n)dx\lo \int_{B_{R_1}} w(x,f(v))f(v)dx.
\end{equation}
Also, it follows from part $ (vii)$  of Proposition 2.1 that the map
$v\rightarrow f(v)$ from $H^1_G$ into $L^q(B_{R} \backslash
B_{R_1})$ is compact for every $q\geq 2,$ hence

\begin{equation}
\int_{B_R\backslash \bar{B}_{R_1}} w(x,f(v_n))f(v_n)dx\lo
\int_{B_R\backslash\bar{B}_{R_1}} w(x,f(v))f(v)dx.
\end{equation}
 Considering (15) and (16), it follows from (14) that
$$\underset{n\rightarrow\infty}{\mathrm{limsup}} \Big |\int w(x,f(v_n))f(v_n)dx-\int
w(x,f(v))f(v)dx \Big | \leq \f{\del}{2},$$ for every $\del>0$.
Consequently
$$\int w(x,f(v_n))f(v_n)dx\lo \int w(x,f(v))f(v)dx,$$
as $n \rightarrow \infty.$  It proves part $(iii).$

To prove  part $(iv),$ note first that $f$ is increasing and $f(0)=0$, hence $f(v_n)\geq 0$ and $ f(v)\geq0.$
For the second term on the right hand side of (8), we have
$$V(x) f(v_n) f'(v_n)\phi \leq  V(x) f(v_n)\phi,$$
and since $v_n\rightharpoonup v$ weakly in $H_1^G$, for the right
hand side of the above inequality we have
$$\lim_{n\rightarrow\infty} \int V(x) f(v_n)\phi \, dx = \int V(x)
f(v)\phi \, dx.$$  Hence by the dominated convergence theorem and
the fact that $v_n \rightarrow v$ a.e. we obtain
\begin{equation}
\lim_{n\rightarrow\infty} \int V(x) f(v_n)f'(v_n)\phi \, dx =  \int V(x) f(v)f'(v)\phi \, dx.
\end{equation}
For the third term on the right hand side of (8), we have
$$w(x,f(v_n)) f'(v_n) \phi \leq \frac{V(x)}{k}f(v_n) \phi, \quad \quad \quad \forall x \in \Lambda^c,$$
and similarly by the dominated convergence theorem, we obtain
 \begin{equation}
\lim_{n\rightarrow\infty} \int_{\Lambda^c} w(x, f(v_n))f'(v_n)\phi \, dx =  \int_{\Lambda^c} w(x, f(v))f'(v)\phi \, dx.
\end{equation}

Also, note that
$$w(x, f(v_n)) f'(v_n) \phi \leq f(v_n)^{22^*-2} \phi   \quad \quad \quad \forall x \in \Lambda,$$
and from part $(vii)$ of Proposition 2.1 that  the map $v\rightarrow
f(v)$ from $H^1_G$ into $L^q(\Lambda)$ is compact for every $q\geq
2,$ hence it follows again from the dominated convergence theorem
that
 \begin{equation}
\lim_{n\rightarrow\infty} \int_{\Lambda} w(x, f(v_n))f'(v_n)\phi \, dx =  \int_{\Lambda} w(x, f(v))f'(v)\phi \, dx.
\end{equation}
It follows from (8) and (17)-(19) that
$$\lim_{n\rightarrow+\infty} \braket{\bar{H}'(v_n)}{\phi}=\braket{\bar{H}'(v)}{\phi}.$$
It proves part $(iv)$. $\square$
\begin{Lem} If $\{v_n\}$ is a $(PS)_{C_0}$ sequence, then $v_n$
converges to $v\in H^1_G$. Consequently
$\bar{H}(v)=\lim_{n\rightarrow+\infty} \bar{H}(v_n)$ and
$\bar{H}'(v)=0$. \end{Lem}
\paragraph{\bf Proof.}  It follows from part $(i)$ of Lemma 3.5 that $v_n$ is a bounded sequence in  $H^1_G.$
Hence, there exists $v \in  H^1_G $ such that, up to a subsequence,
$v_n\rightharpoonup v$ weakly in $H^1_G$ and  $v_n\rightarrow v$
a.e. in ${\BR}^N.$ Since we may replace $v_n$ by $|v_n|,$ we assume
$v_n\geq  0$ and $v\geq 0.$  Since, $\{v_n\}$ is a $(PS)_{C_0}$
sequence we have
\begin{align}
o(\|v_n\|)&= \braket{\bar{H}'(v_n)}{
 \f{f(v_n)}{f'(v_n)}}\\&= \int (1+\f{f(v_n)^2}{1+f(v_n)^2}) |\dDel
v_n|^2dx+ \int V(x)f(v_n)^2dx \quad - \int w(x,f(v_n))f(v_n)dx
\nonumber
\end{align}
and
 \begin{align}
o(\|v\|)=\braket{\bar{H}'(v_n)}{
 \f{f(v)}{f'(v)}}.
\end{align}
It follows from  part $(iv)$ of Lemma 3.5  and (21) that
 \begin{align}
\braket{\bar{H}'(v_n)}{
 \f{f(v)}{f'(v)}} =& \braket{\bar{H}'(v)}{
 \f{f(v)}{f'(v)}}+o(\|v\|) \nonumber \\=&\int (1+\f{f(v)^2}{1+f(v)^2}) |\dDel
v|^2dx+ \int V(x)f(v)^2dx \nonumber \\
& - \int w(x,f(v))f(v)dx +o(\|v\|)
\end{align}
In this step, we show that
\begin{align*}
 \int \frac{f(v)^2 |\dDel v|^2}{1+f(v)^2} \, dx \leq \liminf_{n\rightarrow \infty } \int \frac{f(v_n)^2 |\dDel v_n|^2}{1+f(v_n)^2} \, dx.
\end{align*}
 Set $u_n=f(v_n)$ and $u=f(v).$ A direct
computation shows that
$$\int |\dDel u_n^2|^2 \, dx=4 \int \frac{f(v_n)^2 |\dDel v_n|^2}{1+f(v_n)^2} \, dx\leq 4\|v_n\|^2.$$
Set $w_n=u_n^2.$ It follows from the above that $\{w_n\}_n$ is a
bounded sequence in $D^{1,2}({\BR}^N).$ Hence, up to a subsequence
$w_n \rightharpoonup w$  weakly in $D^{1,2}({\BR}^N)$ and $w_n
\rightarrow w$ a.e. in ${\BR}^N.$ It follows $w=u^2.$ Also, by the
lower semi continuity of the  norm in $D^{1,2}({\BR}^N),$ we obtain
$$\int |\dDel w|^2 \, dx \leq \liminf_{n\rightarrow \infty } \int |\dDel w_n|^2 \, dx.$$
Plug  $w_n=u_n^2$ and $w=u^2$ in this inequality to get
$$\int |\dDel u^2 |^2 \, dx \leq \liminf_{n\rightarrow \infty } \int |\dDel u_n^2 |^2 \, dx.$$
Substituting $u_n=f(v_n)$ and $u=f(v)$ in the above  inequality
gives
\begin{align}
 \int \frac{f(v)^2 |\dDel v|^2}{1+f(v)^2} \, dx \leq \liminf_{n\rightarrow \infty } \int \frac{f(v_n)^2 |\dDel v_n|^2}{1+f(v_n)^2} \, dx.
\end{align}
Also, lower  semi continuity and Fatou's Lemma imply
\begin{align}
\int |\dDel v|^2dx &\leq \underset{n\rightarrow\infty}{\mathrm{liminf}} \int
|\dDel v_n|^2dx, \\
\int V(x)G(v)dx &\leq \underset{n\rightarrow\infty}{\mathrm{liminf}}\int
V(x)G(v_n)dx.
\end{align}
Up to a subsequence  one can assume
\begin{align}
\liminf_{n\rightarrow \infty }\int
|\dDel v_n|^2dx&= \lim_{n\rightarrow \infty } \int
|\dDel v_n|^2dx \\
\liminf_{n\rightarrow \infty }\int
V(x)G(v_n)dx&= \lim_{n\rightarrow \infty }\int
V(x)G(v_n)dx.\\
\liminf_{n\rightarrow \infty } \int \frac{f(v_n)^2 |\dDel v_n|^2}{1+f(v_n)^2} \, dx&=\lim_{n\rightarrow \infty } \int \frac{f(v_n)^2 |\dDel v_n|^2}{1+f(v_n)^2} \, dx.
\end{align}
It follows from (23)-(28) that there exist nonnegative numbers
$\delta_1, \delta_2$ and $\delta_3$ such that
\begin{align}
\lim_{n\rightarrow \infty }\int
|\dDel v_n|^2dx&= \int
|\dDel v|^2dx+\delta_1 \\
\lim_{n\rightarrow \infty }\int
V(x)G(v_n)dx&= \int
V(x)G(v)dx+\delta_2.\\
\lim_{n\rightarrow \infty } \int \frac{f(v_n)^2 |\dDel v_n|^2}{1+f(v_n)^2} \, dx&= \int \frac{f(v)^2 |\dDel v|^2}{1+f(v)^2} \, dx+\delta_3.
\end{align}
Now, we show that $\delta_1=\delta_2=\delta_3=0.$  It follows from
part $(iii)$ of Lemma 3.5 that
$$\int w(x,f(v_n)) f(v_n)dx\lo \int w(x,f(v))f(v)dx.$$
which together with (20) and (22) imply
 \begin{align*}
 \lim_{n\rightarrow \infty } \Big\{\int (1+\f{f(v_n)^2}{1+f(v_n)^2}) |\dDel
v_n|^2dx\\+ \int V(x)f(v_n)^2dx \Big \}&=
\quad   \lim_{n\rightarrow \infty } \int w(x,f(v_n))f(v_n)dx \\
&= \int w(x,f(v))f(v)dx\\
&=\int (1+\f{f(v)^2}{1+f(v)^2}) |\dDel
v|^2dx+ \int V(x)f(v)^2dx
\end{align*}
Taking into account (29), (30) and (31) the above limit implies
$\delta_1=\delta_2=\delta_3=0.$ Therefore, it follows from (29) and
(30) that
\begin{align*}
\int |\dDel v|^2dx&=\underset{n\rightarrow\infty}{\mathrm{lim}} \int
|\dDel
v_n|^2dx \\
\int V(x)G(v)dx&=\underset{n\rightarrow\infty}{\mathrm{lim}}\int
V(x)G(v_n)dx.
\end{align*}
By Proposition 2.1, $v_n\lo v$ in $E_G$ and we have $\dDel v_n\lo
\dDel v$ in $L^2$. Hence $v_n\lo v$ in $H^1_G$. $\square$
\paragraph{\bf Proof of Theorem 3.2.} The proof is a direct
consequence of Lemmas 3.3, 3.4 and 3.5. $\square$
\section{Proof of Theorem 1.1}
To prove Theorem 1.1, note first that every critical point of the
functional  $\bar{J}_\eps$ corresponds to a weak solution of problem
(2). Thus, we need to find a critical point for the functional
$\bar{J}_\eps.$ To do this, we shall show that the functionals
$\bar{J}_\eps$ and $\bar{H}_\eps$ will coincide for the small values
of $\eps$. Hence, every critical point of $\bar{H}_\eps$ will be a
critical point of $\bar{J}_\eps.$ Also, it follows from Theorem 3.2
that $\bar{H}_\eps$ has a
nontrivial critical point for every $\eps>0$. \\

Without loss of generality,  we may assume $\eps^2$ instead of
$\eps$ in the functionals $\bar{H}_\eps$ and $\bar{J}_\eps$, i.e.

$$\bar{H}_\eps (v)=\f{\eps^2}{2} \int |\dDel v|^2+\f {1}{2}\int
V(x)f(v)^2dx-\int W(x,f(v))dx,$$
and

$$\bar {J}_\eps(v)=\f{\eps^2}{2} \int_{{\BR}^N} |\dDel v|^2dx+ \f{1}{2} \int_{{\BR}^N}V(x)
f(v)^2 dx-\f{1}{22^*} \int_{{\BR}^N} |f(v)|^{22^*}dx.$$

It follows from Theorem 3.2 that    there exists a critical point
$v_\eps\in H^G_1$ of $\bar{H}_\eps(v)$ for each $\eps>0$. Set
$u_\eps=f(v_\eps)$.

The following Lemmas are crucial for the proof of Theorem 1.1.

\begin{Lem}
 The sequence $\{u_\eps\}_{\eps>0}$ is strongly convergent
to $0$ when $\eps\lo 0$, in $H^1({\BR}^N)$, i.e.
$$\|u_\eps\|_{H^1}\lo 0 \quad \text{as}\quad \eps\lo 0.$$
\end{Lem}
\paragraph{\bf Proof.} Let $0\not\equiv \phi\in C_{0,r}^\infty({\BR}^N)$ be a
non-negative function with supp$(\phi)\subset\Om$ and $H_1(\phi)\leq
0.$  Set $\ga_1 (t):= h(t \phi).$ Hence, we have $$ \bar { H}_\eps
(\gamma_1 (1))= \bar {H}_\eps (h (\phi))= H_\eps(\phi)\leq
H_1(\phi)\leq0.$$

It follows from the definition of the Mountain Pass value that
$$\bar{H}_\eps(v_\eps)=\inf_{\ga\in\Ga} \sup_{t\in[0,1]}
\bar{H_\eps}(\ga(t))\leq \sup_{t\in[0,1]}
\bar{H_\eps}(\ga_1(t))=\sup_{t\in[0,1]} \bar{H_\eps}(h(t
\phi))=\sup_{t\in[0,1]} {H_\eps}(t\phi).$$

Therefore, we obtain
\begin{align}
\bar{H}_\eps(v_\eps)  &\leq
\sup_{t\in[0,1]}H_\eps(t\phi) \nonumber \\
&=\sup_{t\in [0,1]} \f{\eps^2 t^2}{2} \int |\dDel \phi|^2+\f{\eps^2
t^4}{2} \int
|\phi|^2 |\dDel\phi|^2-\f{t^{22^*}}{22^*} \int |\phi|^{22^*}dx  \nonumber\\
&\leq \sup_{t\in[0,1]}\f{\eps^2 t^2}{2} \int (1+|\phi|^2)
|\dDel\phi|^2dx -
\f{|t|^{22^*}}{22^*}\; \int |\phi|^{22^*}dx \nonumber \\
&\leq (\f{1}{2}-\f{1}{22^*}) \eps^{\f{22^*}{2^*-1}} A(\phi)
\end{align}
where $A(\phi)=\left( \f{\int (1+|\phi|^2) |\dDel \phi|^2dx}{\int
|\phi|^{22^*}dx} \right)^{\f{2^*}{2^*-1}}$. Now, as in the  proof of
part $(i)$ of Lemma 3.5  we obtain
\begin{align}
\bar{H}_\eps(v_\eps) & = \bar{H}_\eps(v_\eps)-\f{1}{\tht} \braket {
\bar{H}'(v_\eps)}{v_\eps} \nonumber \\
&\geq \eps^2(\f{1}{2}-\f{2}{\tht}) \int |\dDel v_n|^2dx+
(\f{1}{2}-\f{1}{\tht})(1-\f{1}{k}) \int V(x)f(v_n)^2dx.
\end{align}
Combining (32) and (33), we get
\begin{equation*}
\eps^2(\f{1}{2}-\f{2}{\tht}) \int |\dDel v_n|^2 dx+
(\f{1}{2}-\f{1}{\tht}) (1-\f{1}{k}) \int V(x) |f(v_n)|^2 dx \leq
(\f{1}{2}-\f{1}{22^*}) \eps^{\f{22^*}{2^*-1}}A(\phi).
\end{equation*}
Therefore
\begin{equation}
(\f{1}{2}-\f{2}{\tht}) \int |\dDel v_n|^2 dx+ (\f{1}{2}-\f{1}{\tht})
(1-\f{1}{k}) \int V(x) f(v_n)^2 dx  \leq (\f{1}{2}-\f{1}{22^*})
\eps^{\f{2}{2^*-1}}A(\phi).
\end{equation}
Hence, substituting $u_\eps=f(v_\eps)$ in (34) implies
$$\int (1+|u_\eps|^2) |\dDel u_\eps|^2dx+\int V(x) |u_\eps|^2dx\leq
C\eps^{\f{2}{2^*-1}}A(\phi).$$
Therefore
$$\|u_\eps\|_{H^1}\lo 0 \quad \text{as}\quad \eps\lo 0.$$
$\square$\\

The following Lemma   is standard (e.g [20]).
\begin{Lem} Let $N>2.$ There is a constant $C=C_N,$ such that
$$|u(x)|\leq \f{C}{|x|^{\f{N-2}{2}}} \|u\|_{H^1({\BR}^N)}\quad
\forall x\neq 0,$$ for any $u \in H_r^1 ({\BR}^N).$
\end{Lem}

\begin{Lem}
 For every compact set $Q\subset{\BR}^N$ such that
$0\not\in Q$, $\|u_\eps\|_{L^\infty(Q)}\lo 0$ as $\eps\lo 0$.
\end{Lem}
\paragraph{\bf Proof.} For each $\eps>0$, it follows from Lemma 4.2 that
$$0\leq u_\eps(x)\leq \f{C}{|x|^{\f{N-2}{2}}} \|u_\eps\|_{H^1({\BR}^N)}\quad
\forall x\neq 0,$$ which together with the result of Lemma 4.1
obviously means
$$\|u_\eps\|_{L^\infty(Q)}\lo 0 \quad \text{as}\quad \eps\lo 0.$$ $\square$
\paragraph{\bf Proof of Theorem 1.1.} By Lemma 4.3 we have
\begin{equation}
M_\eps:=\max_{x\in \bar \La} f(v_\eps)\lo 0 \quad \text{as}\quad
\eps\lo 0.
\end{equation}
%and
%\begin{equation}
%M^2_\eps:=\max_{x\in \pa B_{R_2}} f(v_\eps)\lo 0 \quad
%\textrm{as}\quad \eps\lo 0
%\end{equation}
From (35)  there exists $\eps_0>0$ such that $\max_{x\in \bar
\La}f(v_\eps)<\beta$
%and $\max_{x\in\pa
%B_{R_2}}f_\eps(v_\eps)<\beta$
 for every $0<\eps<\eps_0$. Using the
test function $\phi=\f{(f(v_\eps)-\beta)_+}{f'(v_\eps)}$, we get
\begin{eqnarray*}
0=\braket{\bar{H}'_\eps(v_\eps)}{\phi} = \int_{\textit{F}}
\eps^2(1+\f{f(v_\eps)^2}{1+f(v_\eps)^2}) |\dDel v_\eps|^2 &+&
\int_{{\BR}^N \backslash \bar \La} V(x)f(v_\eps)(f(v_\eps)-\beta)_+dx\\
&-&\int_{{\BR}^N \backslash \bar \La} w(x,f(v_\eps))
(f(v_\eps)-\beta)_+dx
\end{eqnarray*}
where $\textit{F}=({{\BR}}^N \backslash \bar \La)\cap \{x|
f(v_\eps)\geq \beta\}$. From $(g_2)$, we have
$$V(x)f(v_\eps)(f(v_\eps)-\beta)_+- w(x,f(v_\eps))(f(v_\eps)-\beta)_+\geq 0,
\quad \forall x\in\La^c.$$
Thus,
$$\eps^2 \int_{\textit{F}} (1+\f{f(v_\eps)^2}{1+f(v_\eps)^2}) |\dDel v_\eps|^2dx=0,$$
from which we obtain
$$f(v_\eps)\leq \beta, \quad \forall x\in {\BR}^N \backslash \bar \La.$$
Therefore
$$w(x,f(v_\eps))= f(v_\eps)^{22^*-1} , \quad \forall x\in {\BR}^N \backslash \bar \La,$$
and we conclude that
$$\eps^2 \int \dDel v_\eps. \dDel \xi dx+ \int V(x) f(v_\eps) f'(v_\eps) \xi dx =
\int f(v_\eps)^{22^*-1}  f'(v_\eps) \xi dx$$ for every $ \xi\in
H^1_G$ and $\eps\in (0,{\eps}_0)$. Therefore, $\bar {J}_\eps(v)$ has
a critical point $v_\eps$ in $H^1_G$  for every $\eps\in
(0,{\eps}_0)$. $\square$
\begin{Rem} Note that, as in the argument in the proof of Theorem
1.1, it seems the smallness of   $\epsilon$  is required  for
technical reasons. In fact, the smallness of $\epsilon$ ensures that
the deformed functional  $H_{\epsilon}$ and the main functional
$J_{\epsilon}$ coincide and in result they have the same critical
points. However, we don't know if solutions exist for large values
of $\epsilon.$ Indeed, even for the semilinear case ($k=0$),i.e.
\begin{eqnarray}
-\epsilon \Delta u+V(x)u=\mu  |u|^{p-1}u, \quad u>0, x \in
  {\BR}^N,  p+1=2^*,
\end{eqnarray}
 the existence of solutions
depends on the graph topology of coefficient $V(x)$ and the
smallness of $\epsilon.$ In fact, even for this simpler case, it is
not quite clear if solutions exist for large values of $\epsilon.$
\end{Rem}

\paragraph{\textbf{Acknowledgment}:}
 The author is grateful to the referee for his/her critical comments which improved the original
manuscript.\\
%%%%%%%%%%ref.

\end{document}